% sparestyle

\tolerance=10000
%\magnification=1200
\raggedbottom

\baselineskip=15pt
\parskip=1\jot

\def\sk{\vskip 3\jot}

\def\heading#1{\vskip3\jot{\noindent\bf #1}}
\def\label#1{{\noindent\it #1}}

% Reference format

\def\ref#1;#2;#3;#4;#5.{\item{[#1]} #2,#3,{\it #4},#5.}
\def\refinbook#1;#2;#3;#4;#5;#6.{\item{[#1]} #2, #3, #4, {\it #5},#6.} 
\def\refbook#1;#2;#3;#4.{\item{[#1]} #2,{\it #3},#4.}

% Technical notations

\def\({\bigl(}
\def\){\bigr)}

% Greek alphabets

\def\de{\delta}

\def\la{\lambda}

\def\up{\upsilon}

% Calligraphic alphabet

% Boldface alphabet

\def\Stwo#1#2{\left\{{#1\atop#2}\right\}}
\def\Sone#1#2{\left[{#1\atop#2}\right]}

\rightline{\tt cube.arXiv.tex}
\vskip1in

\centerline{\bf The Hypercube of Resistors, Asymptotic Expansions,}
\centerline{\bf  and Preferential Arrangements}
\vskip0.5in

\rightline{NICHOLAS PIPPENGER}
\rightline{Department of Mathematics}
\rightline{Harvey Mudd College}
\rightline{1250 Dartmouth Avenue}
\rightline{Claremont, CA 91711}
\rightline{\tt njp@math.hmc.edu}
\vskip0.5in

\heading{1. Introduction}
\sk

A classic puzzle asks for the resistance between vertices at the ends of a long diagonal when the edges of a cube are replaced by $1$-ohm resistors.
The solution relies on the observation that for each of the endpoints, the three adjacent vertices are at the same potential, by the symmetry of the cube under a $120^\circ$ rotation about the long diagonal.
The network is thus equivalent to one in which three resistors in parallel are in series with six resistors in parallel and with three resistors in parallel, for a total resistance of $1/3 + 1/6 + 1/3 = 5/6$ ohms.
(This problem seems first to have appeared in 1914 in a book by Brooks and Poyser [B4].)
A natural question is:  what happens when the $3$-dimensional cube is replaced by an 
$n$-dimensional hypercube?

Reasoning as before, we observe that all the vertices at a given distance from one of the endpoints
of the long diagonal are again at the same potential, so the network is equivalent to a series connection of parallel connections of resistors.
Since there are $n+1$ distances ($0, 1, \ldots, n$) from one endpoint, there are $n$ parallel connections.
There are ${n\choose k}$ vertices at distance $k$ from the endpoint, and for $0\le k\le n-1$, each of these vertices is connected by $n-k$ resistors to vertices at distance $k+1$.
Thus the total resistance is
$$\eqalignno{
R_n &= \sum_{0\le k\le n-1} {1\over {n\choose k}(n-k)} \cr
& = {1\over n}\sum_{0\le k\le n-1} {1\over{n-1\choose k}}, &(1.1)\cr
}$$
where we have used the identity ${n\choose k}(n-k) = n{n-1\choose k}$, which is easily seen using the expressions for binomial coefficients in terms of factorials.
(This $n$-dimensional version of the problem was posed 
in 1976 by Mullin and Zave as Problem E 2620 in the {\it American Mathematical Monthly\/} [M],
with a solution by Jagers [J],
and again
in 1979 by Singmaster as Problem 79-16 in {\it SIAM Review\/} [S1],
with a solution by Rennie [R].)

The numbers being summed in (1.1) are elements of the ``harmonic triangle'', considered by
Leibniz as a companion to the ``arithmetic triangle'' of Pascal (see Boyer [B3], p.~439).
In the arithmetic triangle,
$$\matrix{
1 \cr
1 &1 \cr
1 &2 &1 \cr
1 &3 &3 &1 \cr
1 &4 &6 &4 &1 \cr
1 &5 &10 &10 &5 &1 \cr
\vdots &&&&&& \ddots ,\cr
}$$
each entry (except the first and  last in each row) is the sum of the elements to its north and its
north-west, whereas in the harmonic triangle,
$$\matrix{
1/1 \cr
1/2 &1/2 \cr
1/3 &1/6 &1/3 \cr
1/4 &1/12 &1/12 &1/4 \cr
1/5 &1/20 &1/30 &1/20 &1/5 \cr
1/6 &1/30 &1/60 &11/60 &1/30 &1/6 \cr
\vdots &&&&&& \ddots ,\cr
}$$
each entry  is the sum of the elements to its south and south-east.
The resistance $R_n$ is the sum of the entries in the $n$-th row:
$$\matrix{
n: &0 &1 &2 &3 &4 &5 &6 &7 &8 &\cdots\cr
R_n: &0 &1 &1 &5/6 &2/3 &8/15 &13/30&151/420 &32/105 &\;\;\cdots. \cr
}$$
The appearance of the reciprocals of binomial coefficients in (1.1) suggests that we also consider
the sum
$$S_n = \sum_{0\le k\le n} {1\over {n\choose k}}, \eqno(1.2)$$
which has the following values:
$$\matrix{
n: &0 &1 &2 &3 &4 &5 &6 &7 &8 &\cdots\cr
S_n: &1 &2 &5/2 &8/3& 8/3 &13/5 &151/60 &256/105 &83/35 &\;\;\cdots. \cr
}$$
Of course, these two sequences are linked by the relations
$$R_n = {1\over n}S_{n-1} \hbox{\ \ \ and\ \ \ } S_n = (n+1)R_{n+1}. \eqno(1.3)$$
In the next section, we shall review some exact results (alternative expressions and generating functions)
for the numbers $R_n$ and $S_n$.
In the following section we shall consider asymptotic expansions for these numbers.
The coefficients in these asymptotic expansions have simple combinatorial interpretations that we shall
consider in the subsequent section.
These combinatorial interpretations will launch us on a tour of old and new results in combinatorial enumeration.
In the final section, we shall return to the hypercube of resistors, and consider the resistance between vertices that are not the endpoints of a long diagonal.
\vfill\eject

\heading{2. Alternative Expressions and Generating Functions}
\sk

The numbers $R_n$ and $S_n$ have alternative expressions,
$$R_n = {1\over 2^n} \sum_{1\le k\le n} {2^k \over k} \eqno(2.1)$$
and
$$S_n = {n+1 \over 2^n} \sum_{0\le k\le n} {2^k \over k+1}, \eqno(2.2)$$
which are equivalent to each other by virtue of (1.3).
While having just as many terms as (1.1) and (1.2), these sums have simpler summands,
and will thus lend themselves more easily to further developments.

The first proof of (2.2) was given by Staver [S3], who derived the recurrence
$S_n = \((n+1)/2n\) S_{n-1} + 1$,
from which (2.2) follows by induction from the base case $S_0 = 1$.
The formula (2.1) was given without proof by Mullin [M], and was given with an ``electrical'' proof by 
Rennie [R], as follows.
First, let a current of $1$ ampere flow into a vertex $A$ and out of a vertex $B$ long-diagonally opposite to $A$.
If $B$ is at potential $0$, then $A$ is at potential $R_n$ volts.
Let $A'$ be adjacent to $A$, and $B'$ long-diagonally opposite to $A'$, and therefore adjacent to $B$.
By symmetry, $1/n$ amperes flows through the $1$-ohm resistor from $A$ to $A'$, so $A'$ is at potential $R_n - 1/n$ volts.
By a similar argument, $B'$ is at potential $1/n$ volts.
Second, reconnect the current source so that $1$ ampere flows into $A'$ and out of $B'$.
If $B'$ is at potential $0$, then $A'$ is at potential $R_n$ volts, $A$ is at potential $R_n - 1/n$ volts,
and $B$ is at potential $1/n$ volts.
Third, suppose that currents of $1$ ampere flow in at each of $A$ and $A'$, and out of each of $B$ and $B'$.
By linearity, we may superimpose the potentials and subtract $1/n$ from their sum,
putting $B$ and $B'$ at potential $0$, and $A$ and $A'$ at potential $2R_n - 2/n$ volts.
But with this final current distribution, $1$ ampere flows from $A$ to $B'$ 
through the resistors of an $(n-1)$-dimensional hypercube, $1$ ampere flows from $A'$ to $B$
through  the resistors of another disjoint  $(n-1)$-dimensional hypercube, and no current flows through
the resistors connecting corresponding vertices in the two hypercubes, since by symmetry they at
at equal potentials.
Thus if $B$ and $B'$ are at potential $0$, $A$ and $A'$ are at potential $R_{n-1}$ volts.
We therefore have $2R_n - 2/n = R_{n-1}$, or $R_n = (1/2)R_{n-1} + 1/n$, from which (2.1) follows by induction from the base case $R_0 = 0$.
Finally, we mention that Sury [S5] proved (2.2) by using the integral representation
$\int_0^1 x^k (1-x)^{n-1-k}\,dx = 1/n{n-1\choose k}$ (Euler's beta integral), summing the resulting
geometric progression inside the integral, and evaluating the resulting integral by a change of variable.

Equations (2.1) and (2.2) allow us to easily derive the generating functions 
$R(z) = \sum_{n\ge 0} R_n\,z^n$ and $S(z) = \sum_{n\ge 0} S_n\,z^n$ for the sequences
$R_n$ and $S_n$.
Indeed, since $-\log(1-z) = z + z^2/2 + z^3/3 + \cdots + z^k/k + \cdots\,$, we see that
$2^k/k$ is the coefficient of $z^k$ in $-\log(1-2z)$.
If $A(z) = \sum_{n\ge 0} A_n\,z^n$ and $B(z) = \sum_{n\ge 0} B_n\,z^n$ are the generating functions for the sequences $A_n$ and $B_n$, respectively, then $C(z) = A(Z)\,B(z)$ is the generating function for
the sequence $C_n = \sum_{0\le k\le n} A_k\,B_{n-k}$, called the ``convolution'' of the sequences
$A_n$ and $B_n$.
As a special case, $B(z) = 1/(1-z)$ is the generating function for the sequence $B_n = 1$, so that
$A(z)/(1-z)$ is the generating function for the sequence $\sum_{0\le k\le n} A_k$ of partial sums of the
sequence $A_n$.
Thus $2^n \sum_{1\le k\le n} 2^k/k$ is the coefficient of $z^n$ in $\(-\log(1-z)\)/(1-z/2)$, so that
$$R(z) = {1\over 1 - z/2}\log{1\over 1-z}. \eqno(2.3)$$
From (1.3), we see that $S(z) = R'(z)$, so differentiating (2.3) yields
$$Sz = {1\over (1-z)(1-z/2)} + {1\over 2(1 - z/2)^2}\log{1\over 1 - z}. \eqno(2.4)$$
Generating functions for sums similar to (2.1) and (2.2) have been given by Pla [P].
\sk

\heading{3. Asymptotic Expansions}
\sk

The results of the preceding section give exact values of $R_n$ and $S_n$  as a rational numbers
but they yield little insight into the behavior of these sequences for large $n$.
To obtain this insight, we develop asymptotic expansions.
It will be convenient to use ``$O$-notation'', where $O\(f(n)\)$ stands for some function $g(n)$
(possibly a different function at each occurrence) such that $\vert g(n)\vert \le c\,f(n)$ for some constant $c$ and all sufficiently large $n$.

We start with (1.1).
Since the binomial coefficients ${n-1\choose k}$ increase as $k$ increases from $0$ to $\lfloor (n-1)/2\rfloor$, then decrease as $k$ increases from $\lceil (n-1)/2\rceil$ to $n-1$, the largest terms in (1.1) are
the first and last: $1/{n-1\choose 0} = 1/{n-1\choose n-1} = 1$. 
The next largest terms are the second and second-to-last, which are 
$1/{n-1\choose 1} = 1/{n-1\choose n-2} = O(1/n)$.
There are $n-4$ other terms, and each of these is at most 
$1/{n-1\choose 2} = 1/{n-1\choose n-3} = O(1/n^2)$,
so the sum of all these other terms is also $O(1/n)$.
Thus we have
$$R_n = {2\over n} + O\left({1\over n^2}\right). \eqno(3.1)$$
This result gives a good estimate of $R_n$ when $n$ is large.

We can refine the estimate (3.1) by extracting the second and second-to-last terms, noting
that the third and third-to-last terms are $O(1/n^2)$, and that each of the remaining $n-6$ terms
is $O(1/n^3)$, so their sum is also $O(1/n^2)$.
This yields
$$R_n = {2\over n}\left(1 + {1\over n-1} + O\left({1\over n^2}\right)\right).$$
Continuing in this way, we obtain
$$R_n = {2\over n}\left(1 + {1\over (n-1)} + {2\over (n-1)(n-2)} + \cdots
 + {k!\over (n-1)(n-2)\cdots(n-k)} + O\left({1\over n^{k+1}}\right)\right) \eqno(3.2)$$
for any fixed $k$.

Equation (3.2) give a sort of asymptotic expansion for $R_n$, but its content would be clearer
if the denominator of each term was a power of $n$, instead of the ``falling powers''  
$(n-1)(n-2)\cdots(n-k)$ that appear there.
That is, we would like an expansion of the form
$$R_n = {2\over n}\left(r_0 + {r_1\over n} + {r_2\over n^2} + \cdots
 + {r_k\over n^k} + O\left({1\over n^{k+1}}\right)\right) \eqno(3.3)$$
 for each $k\ge 0$.
 It is customary to write
 $${n\,R_n\over 2} \sim r_0 + {r_1\over n} + {r_2\over n^2} + \cdots
 + {r_k\over n^k} + \cdots \eqno(3.4)$$ 
 as shorthand for the assertion of (3.3) for each $k\ge 0$.
 The series (3.4) is called an {\it asymptotic expansion}; it is not convergent for any $n$, but it allows
 $R_n$ to be approximated with an error $O(1/n^k)$ for any fixed $k$ and all sufficiently large $n$
 (where the constant hidden in the $O$-notation depends on $k$).
 
 Our task is to determine the coefficients $r_0, r_1,\ldots$ in (3.4).
To do this we shall expand each term $k!/  (n-1)(n-2)\cdots(n-k)$ in (3.2) into a series of 
negative powers of $n$,
$$ {k!\over (n-1)(n-2)\cdots(n-k)}  = \sum_{l\ge k} {k!\,t_{k,l}\over n^l}, \eqno(3.5)$$
then sum the contributions to $r_l$ for each $k\le l$.
First we need to find the numbers $t_{k,l}$ in the expansion (3.5).
These are what have come to be called the ``Stirling numbers of the second kind'', for which we shall
use the notation suggested by Knuth [K, p.~65]: $t_{k,l} = \Stwo{l}{k}$.
These numbers were introduced by James Stirling in 
the Introduction to his {\it Methodus Differentialis\/} [S4] in 1730.
He defined them as the numbers that expand a power $z^l$ of $z$ as a linear combination of the
polynomials $z, z(z-1), \ldots, z(z-1)\cdots(z-l+1)$:
$$z^l = \sum_{0\le k\le l} \Stwo{l}{k}z(z-1)\cdots(z-k+1),$$
and he gave a table for $1\le k\le l\le 9$.
The number $\Stwo{l}{k}$ has a simple combinatorial interpretation: it is the number of ways to
partition the $l$ elements of the set $L = \{1,\ldots, l\}$ into $k$ blocks (non-empty subsets of $L$ that are pairwise disjoint and whose union is $L$).
For $l=3$, for example, we have one partition $\{\{1,2,3\}\}$ into one block,
three partitions $\{\{1\},\{2,3\}\}$, $\{\{1,2\},\{3\}\}$ and $\{\{1,3\},\{2\}\}$ into two blocks and
one partition $\{\{1\},\{2\},\{3\}\}$ into three blocks; thus $\Stwo{3}{1}=1$, $\Stwo{3}{2}=3$ and 
$\Stwo{3}{3}=1$.
At the end of the introduction, Stirling gives the expansion
$${1\over (z+1)(z+2)\cdots(z+k)} = \sum_{l\ge k} (-1)^{l-k}\Stwo{l}{k}{1\over z^l},$$
which, upon substitution of $-n$ for $z$, gives (3.5) in the form
$$ {k!\over (n-1)(n-2)\cdots(n-k)}  = \sum_{l\ge k} k!\Stwo{l}{k}{1\over n^l}.$$
These expansions are in fact convergent for fixed $k\ge 1$ and for $\vert z\vert < 1/k$ or $n > k$
(though Stirling did not distinguish convergent expansions, such as these, and asymptotic expansions, such as (3.4)).
Applying this result to each term in (3.2) gives the desired asymptotic expansion:
$${n\,R_n\over 2} \sim \sum_{l\ge 0} \left(\sum_{0\le k\le l} \Stwo{l}{k}\,k! \right) {1\over n^l}.$$
Thus the coefficients $r_l$ we sought are given by
$$r_l = \sum_{0\le k\le l} \Stwo{l}{k}\,k!. \eqno(3.6)$$

(It will doubtless have occurred to the reader  that if there are ``Stirling numbers of the second kind'', there should also be ``Stirling numbers of the first kind''.
Indeed there are, and they were also introduced by Stirling [S4].
He defined them as the numbers that expand $z(z+1)\cdots(z+l-1)$ as a linear combination of the polynomials $z,z^2, \ldots, z^l$. 
Nowadays it is more common to define them as the 
absolute values of the numbers that expand $z(z-1)\cdots(z-l+1)$ as a linear combination of the polynomials $z,z^2, \ldots, z^l$; in the notation of Knuth [K, p.~65]:
$$z(z-1)\cdots(z-l+1) = \sum_{0\le l\le l} (-1)^{l-k} \Sone{l}{k} z^k.$$
Stirling again gave a table for $1\le k\le l\le 9$, and the expansion
$${1\over z^k} = \sum_{l\ge k} \Sone{l}{k}{1\over z(z+1)\cdots(z+k-1)}.$$
These numbers too have a simple combinatorial interpretation: $\Sone{l}{k}$ is the number of permutations of $l$ elements that have $k$ cycles.
For $l=3$, for example, we have two permutations $(123)$ and $(132)$ with one cycle,
three permutations $(1)(23)$, $(12)(3)$ and $(13)(2)$ with two cycles and
one permutation $(1)(2)(3)$ with three cycles; thus 
$\Sone{3}{1}=2$, $\Sone{3}{2}=3$ and $\Sone{3}{3}=1$.)

 We can find a similar asymptotic expansion for $S_n$.
 Again noting that the largest terms in the sum (1.2) are the first and the last, we obtain
 $$S_n = 2 + O\left({1\over n}\right).$$
 Generalizing this as before yields
 $$S_n = 2\left( 1 + {1\over n} + {2\over n(n-1)} + \cdots + {k!\over n(n-1)\cdots(n-k+1)}
 + O\left({1\over n^{k+1}}\right)\right).$$
Applying (3.5) to each term and summing the contributions for each negative power of $n$, we obtain
$${S_n\over 2}  \sim  s_0 + {s_1\over n} + {s_2\over n^2} + \cdots + {s_k\over n^k}
+ \cdots, \eqno(3.7)$$
where
$$s_l = \sum_{0\le k\le l} \Stwo{l}{k}\,(k+1)!.$$

The coefficients $r_l$ and $s_l$ have simple combinatorial interpretations that we shall study
in the following section.
\sk

\heading{4. Preferential Arrangements}
\sk

In this section we shall study the numbers $r_l$ and $s_l$.
Our model for this study will be a collection of results concerning the ``exponential numbers'' $d_l$,
given by
$$d_l = \sum_{0\le k\le l} \Stwo{l}{k}. \eqno(4.1)$$
These numbers have a simple combinatorial interpretation: $d_l$ is the number of ways to partition
the set $\{1, \ldots, l\}$ into any number of blocks.
For $l=3$, for example, we have seen that there is one partition into one block, three partitions into two blocks and one partition into three blocks; thus $d_3 = 1+3+1 = 5$.
We have the table
$$\matrix{
l: &0 &1 &2 &3 &4 &5 &6 &7 &8 &\cdots\cr
d_l: &1 &1 &2 &5 &15 &52 &203 &877 &4140 &\;\;\cdots. \cr
}$$
(The sequence $d_l$ is A000110 in Sloan's 
{\it On-Line Handbook of integer Sequences\/}  [S2].)

There are three aspects of the exponential numbers that are of particular interest to us:
a recurrence, a generating function and an expression as the sum of an infinite series.
The recurrence is
$$d_l = \de_l + \sum_{0\le k\le l-1} {l-1\choose k} d_k, \eqno(4.2)$$
where $\de_l$ is $1$ for $l=0$ and $0$ for all other values of $l$.
This recurrence allows $d_l$ to be computed from the previous values $d_0, d_1, \ldots, d_{l-1}$.

The ``exponential generating function'' $d(z) = \sum_{l\ge 0} d_l \, z^l / l!$ is given by
$$d(z) = e^{e^z - 1}, \eqno(4.3)$$
where the term ``exponential'' refers to the factor $1/l!$ in the defining sum (in contrast to the 
``ordinary'' generating functions that we used in Section 2).

The expression as an infinite sum is
$$d_l ={1\over e}  \sum_{n\ge 0} {n^l \over n!}, \eqno(4.4)$$
where $e = 2.7182\ldots$ is the base of natural logarithms.
Since $d_l$ is expressed as a finite sum in (4.1), it may not be clear what advantage there is to (4.4).
But the finite sum involves the Stirling numbers of the second kind, whereas the infinite sum involves only powers and factorials.
Furthermore, the sum (4.4) has a simple probabilistic interpretation:
$d_l$ is the $l$-th moment ${\rm Ex}[N^l]$ of a Poisson-distributed random variable $N$ with mean 
$\la=1$ (since $\Pr[N=n] = e^{-\la} \, \la^n / n! = 1/en!$ for such a random variable).

The exponential numbers were mentioned in 1934 by Bell [B1, B2], and are on that account sometimes called the ``Bell numbers''.
But (4.2), (4.3) and (4.4) were all given earlier:
(4.2) in 1933 by Touchard [T] (in an equivalent ``umbral'' form), (4.3) in 1886 by Whitworth 
[W, Proposition XXIV, p.~95] and (4.4) in 1877 by Dobi\'{n}ski [D] (who actually only gave the cases
$1\le l\le 8$; but it is clear from his derivations that $d_l$ satisfies the recurrence (4.2)).

In his marvelous book {\it Asymptotic Methods in Analysis}, N.~G. de Bruijn [B5, Section 3.3]
derives the asymptotic expansion 
$${1\over n!} \sum_{0\le k\le n} k! \sim  1 + {d_0\over n^1} + {d_1\over n^2} + \cdots + 
{d_k \over n^{k+1}} + \cdots,$$
and then says that it is ``only for the sake of curiosity'' that he mentions that the coefficients $d_l$,
given by (4.1),
have a combinatorial interpretation.
One of our goals in this paper is to pursue this curiosity; our motto is: whenever the coefficients in 
an expansion are integers, look for a combinatorial interpretation!

The coefficients $r_l$, given by (3.6),
also have a simple combinatorial interpretation: they are the number of ways of ranking $l$ candidates,
with ties allowed; that is, the $l$ candidates are first to be partitioned into equivalence classes,
then the equivalence classes are to be linearly ordered.
This interpretation follows from those of $\Stwo{l}{k}$ and $k!$, where $k$ is the number of equivalence classes in the partition.
Because of this interpretation, $r_l$ is called the number of {\it preferential arrangements\/} of $l$ elements.
For $l=3$, for example, the one partition into one block can have its block ordered in one way,
each of the three partitions into two blocks can have its blocks ordered in two ways and the partition into three blocks can have its blocks ordered in six ways; thus $r_3 = 1\cdot 1 + 3\cdot 2 + 1\cdot 6 = 13$.
We have the table
 $$\matrix{
l: &0 &1 &2 &3 &4 &5 &6 &7 &8 &\cdots\cr
r_l: &1 &1 &3 &13 &75 &541 &4683 &47293 &545835 &\;\;\cdots. \cr
}$$
(The sequence $r_l$ is A000670 in Sloan [S2].)
We shall derive the recurrence
$$r_l = \de_l + \sum_{0\le k\le l-1} {l\choose k} r_k, \eqno(4.5)$$
the exponential generating function (defined by  $r(z) = \sum_{l\ge 0} r_l \, z^l / l!$)
$$r(z) = {1\over 2 - e^z} \eqno(4.6)$$
and the summation expression
$$r_l = {1\over 2} \sum_{n\ge 0} {n^l \over 2^n}. \eqno(4.7)$$

We begin by deriving the recurrence (4.5).
For $l\ge 1$, we can construct a preferential arrangement on $l$ candidates by first choosing
the number $k$ of candidates tied in the top equivalence class (with $k$ in the range $1\le k\le l$), then 
choosing (in one of ${l\choose k}$ ways) the  candidates in this class, and finally choosing 
(in one of $r_{l-k}$ ways) a preferential arrangement of the remaining $l-k$ candidates.
This gives the recurrence
$$\eqalignno{
r_l 
&= \sum_{1\le k\le l} {l\choose k}\,r_{l-k} \cr
&= \sum_{0\le j\le l-1} {l\choose j+1}\,r_{l-j-1} \cr
&= \sum_{0\le k\le l-1} {l\choose l-k}\,r_k, \cr
&= \sum_{0\le k\le l-1} {l\choose k}\,r_k, \cr
}$$
where we first made the substitution $k = j+1$, then the substitution $j = l-1-k$, and finally used the identity ${l\choose l-k} = {l\choose k}$.
This equation  holds for $l\ge 1$; since $r_0 = 1$, we obtain (4.5) for $l\ge 0$.

Next we shall derive the
exponential generating function (4.6).
Adding $r_l$ to both sides of (4.5) yields
$$2r_0 =  \de_l + \sum_{0\le j\le l} {l\choose j}\,r_j. \eqno(4.8)$$
Multiplying both sides of this equation by $z^l / l!$ and summing over $l\ge 0$, we obtain
$$\eqalignno{
2r(z)
&= 1 + \sum_{l\ge 0} {z^l\over l!} \sum_{0\le j\le l} {l\choose j}\,r_j \cr
&= 1 + \sum_{j\ge 0} {z^j \, r_j \over j!} \sum_{k\ge 0} {z^k\over k!} \cr
&= 1 + e^z \, r(z), &(4.9)\cr
}$$
where we have made the substitution $k = l-j$ and used the identity $e^z = \sum_{k\ge 0} z^k / k!$.
Solving this equation for $r(z)$ yields (4.6).

Finally we shall derive the summation expression (4.7).
To do this, we rewrite the exponential generating function $r(z)$ from (4.6):
$$\eqalign{
\sum_{l\ge 0} {r_l \, z^l \over l!}
&= {1\over 2} \, {1\over 1 - {1\over 2}e^z} \cr
&= {1\over 2} \sum_{n\ge 0} {e^{nz} \over 2^n} \cr
&= {1\over 2} \sum_{n\ge 0} {1 \over 2^n} \sum_{l\ge 0} {(nz)^l \over l!} \cr
& = \sum_{l\ge 0} \left({1\over 2} \sum_{n\ge 0} {n^l \over 2^n}\right) {z^l \over l!}. 
}$$
Since the coefficient of $z^l / l!$ must be the same on both sides of this equation, we obtain (4.7).

The name ``preferential arrangement'' was introduced by Gross [G], as was the summation expression (4.7).
The numbers $r_l$ (with a different combinatorial interpretation involving trees), the recurrence (4.5) and the generating function (4.6) were given by Cayley [C1] in 1859; the combinatorial interpretation
we have used is implicit in 1866 by Whitworth [W, Proposition XXII, p.~93] (Whitworth shows that the terms $\Stwo{l}{k}\,k!$ for fixed $k$ have the exponential generating function $(e^z - 1)^k$; summation over $k\ge 0$ then yields $1/\(1 - (e^z - 1)\) = 1/(2 - e^z)$.)

We turn now to the numbers $s_l$, which also have a simple combinatorial interpretation:
$s_l$ is the number of ways of ranking $l$ candidates, with ties allowed, and with a ``bar''
that may be placed above all the candidates, between two equivalence classes of tied candidates, or below all the candidates.
Thus we may call $s_l$ the number of {\it barred preferential arrangements\/} of $l$ elements.
If there are $k$ equivalence classes of tied candidates, there are $k+1$ positions for the bar.
For $l=3$, for example, the one preferential arrangement with one block has two positions for the bar,
each of the six preferential arrangements with two blocks has three positions for the bar, and each of the six preferential arrangements with three blocks has four positions for the bar; thus
$s_3 = 1\cdot 2 + 6\cdot 3 + 6\cdot 4 = 44$.
We have the table
$$\matrix{
l: &0 &1 &2 &3 &4 &5 &6 &7 &8 &\cdots\cr
s_l: &1 &2 &8 &44 &308 &2612 &25988 &296564 &3816548 &\;\;\cdots. \cr
}$$
(The sequence $s_l$ is A005649 in Sloan [S2].)

In lieu of a recurrence for the numbers $s_n$, we shall derive a formula expressing them in terms of the
numbers $r_l$:
$$\eqalignno{
s_l
&= \sum_{0\le k\le l} {l\choose k} r_k \, r_{l-k}. &(4.10)
}$$
We shall also derive the exponential generating function 
(defined by  $s(z) = \sum_{l\ge 0} s_l \, z^l / l!$)
$$s(z) = {1\over (2 - e^z)^2} \eqno(4.11)$$
and the summation expression
$$s_l = {1\over 4} \sum_{n\ge 0} {(n+1)n^l \over 2^n}. \eqno(4.12)$$

For $l\ge 0$, we can construct a barred preferential arrangement on $l$ candidates by first choosing
the number $k$ of candidates above the bar (with $k$ in the range $0\le k\le l$), then 
choosing (in one of ${l\choose k}$ ways) the  candidates above the bar, then choosing
(in one of $r_k$ ways) a preferential arrangement of these candidates,
and finally choosing (in one of $r_{l-k}$ ways)
a preferential arrangement of the remaining $l-k$ candidates.
This gives the formula (4.10).
Next, multiplying both sides of (4.10) by $z^l / l!$ and summing over $l\ge 0$ yields
$$\eqalignno{
s(z)
&= \sum_{l\ge 0} {z^l \over l!}  \sum_{0\le k\le l} {l\choose k} r_k \, r_{l-k} \cr
&= \sum_{k\ge 0} {z^k \, r_k \over k!} \sum_{j\ge 0} {z^j \, r_j \over j!} \cr
&= r(z)^2, &(4.13)\cr
}$$
where we have made the substitution $k = l-j$.
Substituting (4.6) in this equation yields (4.11).
Finally, reasoning similar to that used to derive (4.7) leads to (4.12).

Before concluding this section, let us derive two more identities relating $r_l$ and $s_l$:
$$r_{l+1} = \sum_{0\le k\le l} {l\choose k} \, s_k \eqno(4.14)$$
and
$$r_l + r_{l+1} = 2s_l. \eqno(4.15)$$
These can be given direct combinatorial proofs (and the reader may enjoy finding these), but we shall use two different methods that are often useful when dealing with sequences that have explicit exponential generating functions.

To prove (4.14), we use the notion of ``binomial convolution''.
Suppose that $a(z) = \sum_{l\ge 0} a_l\,z^l/l!$ and $b(z) = \sum_{l\ge 0} b_l\,z^l/l!$ are the exponential generating functions for the sequences $a_l$ and $b_l$, respectively.
Then 
$$\eqalign{
a(z)\,b(z) 
&= \sum_{k\ge 0} {a_k\,z^k\over k!} \; \sum_{j\ge 0} {b_j\,z^j\over j!} \cr
&= \sum_{l\ge 0} {z^l\over l!} \; \sum_{0\le k\le l} {l\choose k}\,a_k\,b_{l-k}, \cr
}$$
where we have made the substitution $j=l-k$.
Thus $c(z) = a(z)\,b(z)$ is the exponential generating function for the sequence 
$c_l = \sum_{0\le j\le l} {l\choose j}\,a_j\,b_{l-j}$, which is called the {\it binomial convolution\/}
of the sequences $a_l$ and $b_l$, and denoted $(a*b)_l$.
(We have already encountered binomial convolutions twice in this section: once to derive
(4.9) from (4.8), where the convolution can be expressed as $s = \de + r*\up$ (the sequence
$\up_l = 1$ for all $l\ge 0$, and has exponential generating function $\up(z) = e^z$), and again to derive (4.13) from (4.10).)
We shall also need the fact that $a'(z)$ (where the prime indicates differentiation) is the exponential generating function for the sequence $a_{l+1}$, which we shall denote $a'_l$.

To derive (4.14), we may now observe that $r'(z) = e^z / (2-e^z)^2 = e^z \, s(z)$.
Thus $r_{l+1} = r'_l = (s*\up)_l$, which yields (4.14).

To derive (4.15), we note that $r'(z) = e^z / (2-e^z)^2$ implies that
$r(z)$ satisfies the differential equation
$$r'(z) + r(z) = 2r(z)^2. \eqno(4.16)$$
(This differential equation, together with the initial condition $r(0) = 1$, uniquely determines $r(z)$.
In fact, it is an example of  a ``Riccati equation'', which can be solved to analytically for $r(z)$.)
Substituting (4.13) in (4.16), we obtain (4.15).
(We note that (4.15) can also be obtained from (4.7) and (4.12).)

That the numbers $s_l$, defined by (3.7), have the exponential generating function given in (4.11) was given as an exercise (without proof or reference) by Comtet [C2, p.~294, Ex.~15].
Our combinatorial interpretation of these numbers in terms of barred preferential arrangements seems to be new.
\sk

\heading{5. More Resistances}

We mentioned in the introduction that Singmaster [S1] posed in 1978 the problem of determining $R_n$.
What we did not mention then is that he asked not only for $R_n$, but for $R_{n,k}$, the resistance between two vertices at distance $k$ (for $1\le k\le n$) in an $n$-dimensional hypercube of $1$-ohm resistors.
Rennie's solution [R] to Singmaster's problem covered 
(by various arguments) the cases $k = 1, 2$ and $k = n, n-1, n-2$.
In this section we shall give (by a single argument) a complete solution to Singmaster's problem: 
for $0\le k\le n$,
$$R_{n,k} = 
{2\over n} \sum_{0\le j\le k-1} {1\over {n-1\choose j}} \;  
{1\over 2^n} \sum_{j+1\le i\le n} {n\choose i}. \eqno(5.1)$$
Our solution, like that of Rennie, is based on the principles of symmetry and superposition.

Consider the situation in which a current of $1$ ampere flows out of a vertex $A$, while currents of
$1/(2^n - 1)$ amperes flow into each of the $2^n - 1$ other vertices.
Symmetry ensures that all ${n\choose j}$ vertices at distance $j$ from $A$ are at the same potential.
Call this potential $U_j$ volts, where $U_0 = 0$.
There are $1/n{n-1\choose j}$ $1$-ohm resistors connecting vertices at potential $U_j$ to vertices at potential $U_{j+1}$, and a total current of $\sum_{j+1\le i\le n} {n\choose i}/(2^n - 1)$ amperes flows through them.
By Ohm's law,
$$U_{j+1} - U_j = {1\over n{n-1\choose j}} \; {1\over 2^n - 1} \sum_{j+1\le i\le n} {n\choose i}, $$
and thus
$$U_k = {1\over n} \sum_{0\le j\le k-1}
{1\over {n-1\choose j}} \; {1\over 2^n - 1} \sum_{j+1\le i\le n} {n\choose i}, $$
Now let $B$ be a vertex at distance $k$ from $A$, and consider the situation in which 
a current of $1$ ampere flows into $B$ and currents of $1/(2^n - 1)$ amperes flow out of  each of the $2^n - 1$ other vertices.
In this situation there is again a potential difference of $U_k$ volts between $A$ and $B$.
By superposition, if a current of $1 + 1/(2^n - 1)$ amperes flows into $B$ and out of $A$,
there will be a potential difference of $2U_k$ between these vertices.
Again using  Ohm's law, we have
$$R_{n,k} = {2\over 1 + 1/(2^n - 1)}\,
{1\over n} \sum_{0\le j\le k-1}
{1\over {n-1\choose j}} \; {1\over 2^n - 1} \sum_{j+1\le i\le n} {n\choose i},$$
which yields (5.1).

Let us first check that (5.1) agrees with (1.1) for $k=n$.
We have
$$\eqalign{
R_{n,n}
&= {2\over n} \sum_{0\le j\le n-1} {1\over {n-1\choose j}} \;  
{1\over 2^n} \sum_{j+1\le i\le n} {n\choose i} \cr
&= {1\over n} \sum_{0\le j\le n-1} {1\over {n-1\choose j}} \;  
{1\over 2^n} \sum_{j+1\le i\le n} {n\choose i} \cr
&\qquad + {1\over n} \sum_{0\le j\le n-1} {1\over {n-1\choose j}} \;  
{1\over 2^n} \sum_{j+1\le i\le n} {n\choose i} \cr
&= {1\over n} \sum_{0\le j\le n-1} {1\over {n-1\choose j}} \;  
{1\over 2^n} \sum_{j+1\le i\le n} {n\choose i} \cr
&\qquad + {1\over n} \sum_{0\le j\le n-1} {1\over {n-1\choose j}} \;  
{1\over 2^n} \sum_{0\le i\le j} {n\choose i} \cr
&= {1\over n} \sum_{0\le j\le n-1} {1\over {n-1\choose j}} \;  
{1\over 2^n} \sum_{0\le i\le n} {n\choose i} \cr
&= {1\over n} \sum_{0\le j\le n-1} {1\over {n-1\choose j}}  \cr
&= R_n, \cr
}$$
where we have used the identities ${n-1\choose j} = {n-1\choose n-1-j}$ and 
${n\choose i} = {n\choose n-i}$.

Next, let us show that $R_{n,k}$ is an increasing function of $k$ (that is, that 
$\nabla_k R_{n,k} = R_{n,k} - R_{n,k-1} > 0$ for $1\le k\le n$).
From (5.1) we have
$$\nabla_k R_{n,k} = 
 {1\over n 2^{n-1}} \;{1\over {n-1\choose k-1}} \;   \sum_{k\le i\le n} {n\choose i}, \eqno(5.2)$$
 and the expression on the right-hand side is obviously positive.

Finally, let shall show that $R_{n,k}$ is a concave function of $k$ (that is, that 
$\nabla_k^2 R_{n,k} = R_{n,k} - 2R_{n,k-1}  + R_{n,k-2} < 0$ for $2\le k\le n$).
From (5.2) we have
$$\nabla_k^2 R_{n,k} = 
 {1\over n 2^{n-1}} \left({1\over {n-1\choose k-1}} \;   \sum_{k\le i\le n} {n\choose i}
-  {1\over {n-1\choose k-2}} \;  \sum_{k-1\le i\le n} {n\choose i}\right).$$
The expression on the right-hand side is obviously negative for 
$k \le n/2 + 1$, since in this case we have
 ${n-1\choose k-2} \le  {n-1\choose k-1}$ and 
 $\sum_{k\le i\le n} {n\choose i} <  \sum_{k-1\le i\le n} {n\choose i}$.
 For the case $k > n/2 + 1$, we factor  ${n-1\choose k-2}$ out of the expression in parentheses, then move all but the first term from the second sum into the first sum:
$$\eqalign{
\nabla_k^2 R_{n,k} 
&=  {1\over n2^{n-1}} \; {1\over {n-1\choose k-2}}
\left({k-1\over n-k+1} \;  \sum_{k\le i\le n} {n\choose i} - 
 \sum_{k-1\le i\le n} {n\choose i}\right) \cr
 &=  {1\over n2^{n-1}} \; {1\over {n-1\choose k-2}}
 \left({2k - n - 2\over n-k+1} \;  \sum_{k\le i\le n} {n\choose i} - 
 {n\choose k-1}\right). \cr
 }$$ 
Thus what remains to be proved is that the expression in parenthesis is negative; that is, that
$${2k - n - 2\over n-k+1} \;  \sum_{k\le i\le n} {n\choose i} <
 {n\choose k-1}$$
for $k > n/2 + 1$,
 or equivalently, by the substitution $k=n-j$ and the identity ${n\choose i} = {n\choose n-i}$,
$${n-2j-2 \over j+1} \sum_{0\le i\le j} {n\choose i} < {n\choose j+1} \eqno(5.3)$$
for $0\le j < n/2 - 1$.

To prove (5.3), we observe that $0\le i < n/2 - 1$ implies that
$${n\choose i} \le {i+1 \over n-i} {n\choose i+1}.$$
Since $(i+1)/(n-i)$ is an increasing function of $i$, we have
$${n\choose i} \le \left({j+1 \over n-j}\right)^{j-i+1} {n\choose j+1}$$
for $i\le j < n/2 - 1$.
Thus we may bound the sum in (5.3) by the sum of an infinite geometric series,
$$\eqalign{
 \sum_{0\le i\le j} {n\choose i}
 &\le  \sum_{0\le i\le j}  \left({j+1 \over n-j}\right)^{j-i+1} {n\choose j+1} \cr
 &\le  \sum_{m\ge 1}  \left({j+1 \over n-j}\right)^{m} {n\choose j+1} \cr
 &= {j+1 \over 2 -2j - 1} {n\choose j+1}. \cr
 }$$
 This inequality proves (5.3), and thus completes the proof that $R_{n,k}$ is concave.
 \sk
 
 \heading{6. References}
 \sk
 
 \ref B1; E. T. Bell;
 ``Exponential Polynomials'';
 Annals of Mathematics; 35:2 (1934) 258--277.
 
 \ref B2; E. T. Bell;
 `Exponential Numbers'';
 American Mathematical Monthly; 41:7 (1934) 411--419.
 
 \refbook B3; C. B. Boyer;
 A History of Mathematics;
 John Wiley \& Sons, 1968.
 
 \refbook B4; E. E. Brooks and A. W. Poyser;
 Electricity and Magnetism:
 A Manual for Advanced Classes;
 Longmans, Green, 1914.
 % p. 277, 5/6 and 7/12
 
 \refbook B5; N. G. de Bruijn;
 Asymptotic Methods in Analysis;
 North-Holland, 1958 (reprinted by Dover, 1981).
 % Section 3.3
 % There is usually no reason to try to obtain an explicit for the coefficients of a divergent
 % asymptotic series. ... So it is only for the sake of curiosity that we mention that ....
 
 \ref C1; A. Cayley;
 ``On the Analytical Forms Called Trees---Part II'';
 Philosophical Magazine; 18 (1859) 374--378.
 
 \refbook C2; L. Comtet;
 Advanced Combinatorics;
 D.~Reidel Publishing, 1974.
 % p. 228, Ex. 20
 % a_m = \sum_{1\le k\le m} k! S2(m,k), \sum_{m\ge 0} a_m z^m/m! = 1/(2-e^z)
 %
 % p. 294, Ex. 15
 %  I_n = \sum_{0\le k\le n} {n\choose k}^{-1} = (n+1)2^{-(n+1)}\sum_{1\le k\le n+1} 2^k/k
 % 2I_n = ((n+1)/2n)I_{n-1} + 2
 % I_n/2 \sim 1 + \sum_{p\eg 0} b_p n^{-p-1}
 % \sum_{p\ge 0} b_p z^p/p! = 1/(2-e^z)^2
 
 \ref D; G. Dobi\'{n}sky;
 ``Summierung der Reihe $\sum n^m / n!$ f\"{u}r $m=1,2,3,4,5,\ldots$'';
 Archiv f\"{u}r Mathematik und Physik ($=$ Grunert's Archiv); 61 (1877) 333--336.

 \ref G; O. A. Gross;
 ``Preferential Arrangements'';
 American Mathematical Monthly; 69:1 (1962) 2--4. 
 
 \refbook I; E. L. Ince;
 Ordinary Differential Equations;
 Dover, 1956.
 % Ricatti: pp. 23--25.
 
 \ref J; A. A. Jagers;
 ``Solution to Elementary Problem E 2620'' (also solved by T.~Morley and D.~Zave);
 American Mathematical Monthly; 85:2 (1978) 117--118.
 
  \refbook K; D. E. Knuth;
 The Art of Computer Programming, v.~1: Fundamental Algorithms;
 Addison-Wesley, 1968.
 
 \ref M; A. Mullin and D. Zave;
 ``Elementary Problem E 2620'';
 American Mathematical Monthly; 83 (1976) 740.

 \ref P; J. Pla;
 ``The Sum of the Inverses of the Binomial Coefficients Revisited'';
 Fibonacci Quarterly; 35: 4 (1997) 342--345.
 
 \ref R; B. C. Rennie;
 ``Solution to Problem 79-16'';
 SIAM Review; 22 (1980) 504-508.
  
 \ref S1; D. Singmaster;
 ``Problem 79-16'';
 SIAM Review; 21:4 (1978) 559.
 
 \item{[S2]} N. J. A.  Sloane (Ed.),
 {\it The On-Line Encyclopedia of Integer Sequences},
 published electronically at 
 {\tt http://www.research.att.com/\~{}njas/sequences}.
 % r_n A000670
 % s_n A005649
 
 \ref S3; T. B. Staver;
 ``Om summasjon av potenser av binomiaalkoeffisientene'';
 Norsk Matematisk Tidsskrift; 29 (1947) 97--103.
 
 \refbook S4; J. Stirling;
 Methodus Differentialis;
 1730 (Translated and annotated by I.~Tweddle, Springer, 2003).
 
 \ref S5; B. Sury;
 ``Sum of the Reciprocals of the Binomial Coefficients'';
 European J. Combinatorics; 14:4 (1993) 351--353.
 
 \ref T; J. Touchard;
 ``Propri\'{e}t\'{e}s arithm\'{e}tiques de certains nombres r\'{e}currents'';
 Ann.\ Soc.\ Sci.\ Bruxelles; 53 (1933) 21--31.
 
\refbook W; W. A. Whitworth;
Choice and Chance;
Deighton, Bell and Co., 1886.

\bye